\documentclass[twocolumn,10pt]{asme2e_m}

\usepackage[cmex10]{amsmath}
\usepackage{verbatim}
\usepackage[table]{xcolor}
\usepackage{tikz}
\usepackage{gensymb}
\usepackage{soul}
\usepackage{graphicx}
\usepackage{amsmath}
\usepackage{array}
\usepackage{amssymb}
\usepackage{verbatim}
\usepackage{cleveref}
\usepackage{color}
\usepackage{relsize}
\usepackage{tikz}
\usepackage{cite}
\usepackage{epsfig}
\usepackage{fancyhdr}
\usepackage{graphics}
\usepackage{float}
\usepackage{tablefootnote}
\confshortname{DSCC 2016} \conffullname{ASME 2016 Dynamic
Systems and Control Conference} \confdate{12-14}
\confmonth{October} \confyear{2016} \confcity{Minneapolis, Minnesota} \confcountry{USA}
\papernum{DSCC2016-9732}
\title{Handling Model and Implementation Uncertainties \\ via An Adaptive Discrete Sliding Mode Controller Design}
\author{Mohammad Reza Amini\thanks{Address all correspondence to this author.} \\
       {\tensfb Mahdi Shahbakhti}
       \affiliation{
Dept. of Mechanical Engineering\\
Michigan Technological University\\
Houghton, MI 49931\\
Emails: mamini, mahdish@mtu.edu
    }
}
\author{Selina Pan
    \affiliation{Dept. of Mechanical Engineering\\
Stanford University\\
Stanford, CA 94305\\
Email: slpan@stanford.edu
    }
}
\author{J. Karl Hedrick
    \affiliation{Dept. of Mechanical Engineering\\
University of California\\
Berkeley, CA 94720\\
Email: khedrick@me.berkeley.edu
    }
}
\graphicspath{{./figures/}}
\begin{document}
\maketitle
\thispagestyle{empty}
\pagestyle{empty}
%
%%%%%%%%%%%%%%%%%%%%%%%%%%%%%%%%%%%%%%%%%%%%%%%%%%%%%%%%%%%%%%%%%%%%%%
\begin{abstract}
%There are two main sources of uncertainties in design cycle of industrial model-based controllers. First, the imprecisions arise during digital implementation of the controllers upon analog-to-digital (ADC) conversion of signals at the input/output (I/O) of the controller due to sampling and equalization. Second, the uncertainty in the modeled plant's dynamics used in the controller formulation. 
Analog-to-digital conversion (ADC) and uncertainties in modeling the plant dynamics are the main sources of imprecisions in the design cycle of model-based controllers. These implementation and model uncertainties should be addressed in the early stages of the controller design, otherwise they could lead to failure in the controller performance and consequently increase the time and cost required for completing the controller verification and validation (V\&V) with more iterative loops. In this paper, a new control approach is developed based on a nonlinear discrete sliding mode controller (DSMC) formulation to mitigate the ADC imprecisions and model uncertainties. To this end, a DSMC design is developed against implementation imprecisions by incorporating the knowledge of ADC uncertainties on control inputs via an online uncertainty prediction and propagation mechanism. Next, a generic online adaptive law will be derived to compensate for the impact of an unknown parameter in the controller equations that is assumed to represent the model uncertainty.
%in the modeled dynamics of the plant. 
The final proposed controller is an integrated adaptive DSMC with robustness to implementation and model uncertainties that includes (i) an online ADC uncertainty mechanism, and (ii) an online adaptation law. The proposed adaptive control approach is evaluated on a nonlinear automotive engine control problem in real-time using a processor-in-the-loop (PIL) setup with an actual electronic control unit (ECU). The results
reveal that the proposed adaptive control technique removes the uncertainty in the model fast, and significantly improves the robustness of the controllers to ADC imprecisions. This
provides up to 60\% improvement in the performance of the  controller under implementation and model uncertainties compared to a baseline DSMC, in which there are no incorporated ADC imprecisions. 
\end{abstract}

%%%%%%%%%%%%%%%%%%%%%%%%%%%%%%%%%%%%%%%%%%%%%%%%%%%%%%%%%%%%%%%%%%%%%%
\begin{nomenclature}\hspace{-8.5mm}% \vspace{-0.5cm}
\begin{tabular*}{0.95\columnwidth}{ l p{8cm}}
  \multicolumn{2}{l}{   } \\
$AFR$ &air fuel ratio, [-]\\ %\vspace{-0.16cm}
%$Cp$ &heat capacity at constant pressure, [$kJ/(kg.K)$]\\ %\vspace{-0.16cm}
%$\dot{HC}_{eng}$ &rate of engine raw HC emission, [$kg/s$]\\ %\vspace{-0.16cm}
%$\dot{HC}_{tp}$ &rate of tailpipe HC emission, [$kg/s$]\\ %\vspace{-0.16cm}
$J$ &effective engine inertia, [$kg.m^2$]\\ %\vspace{-0.16cm}
%$m$ &mass of the catalytic converter, [$kg$]\\ %\vspace{-0.16cm}
$m_a$ &air mass inside the intake manifold, [$kg$]\\ %\vspace{-0.16cm}
$\dot{m}_f$ &fuel mass flow into the cylinder, [$kg/s$]\\ %\vspace{-0.16cm}
$\dot{m}_{ai}$ &air mass flow into the intake manifold, [$kg/s$]\\ %\vspace{-0.16cm}
$\dot{m}_{ao}$ &air mass flow into the cylinders, [$kg/s$]\\ %\vspace{-0.16cm}
$\dot{m}_{fc}$ &commanded fuel mass flow, [$kg/s$]\\ %\vspace{-0.16cm}
%$\dot{Q}_{in}$ &heat flow into the catalytic converter, [$J/s$]\\ %\vspace{-0.16cm}
%$\dot{Q}_{out}$ &heat flow from the catalytic converter, [$J/s$]\\ %\vspace{-0.16cm}
%$\dot{Q}_{gen}$ &generated heat flow in the cat. converter, [$J/s$]\\ %\vspace{-0.16cm}
%$r_c$ &compression ratio, [-]\\ %\vspace{-0.16cm}
%$T_{cat}$ &catalytic converter temperature, [$K$]\\ %\vspace{-0.16cm}
$T_{exh}$ &exhaust gas temperature, [$K$]\\ %\vspace{-0.16cm}
%$T_E$ &Engine torque, [$N.m$]\\\vspace{-0.17cm}
$\Delta$ &spark timing, [$^{o}$~$ATDC$]\\ %\vspace{-0.16cm}
%$\eta_{cat}$ &catalytic converter conversion efficiency, [-]\\ %\vspace{-0.16cm}
%$\eta_{vol}$  &engine volumetric efficiency, [-]\\ %\vspace{-0.16cm}
%$\lambda$ &gradient of sliding surface, [-]\\ %\vspace{-0.16cm}
$\omega_e$  &engine speed, [$rad/s$] \\ %\vspace{-0.16cm}
$N$  &engine speed, [$RPM$] \\ %\vspace{-0.16cm}
$\tau_{f}$ &fuel evaporation time constant, [$s$]\\ %\vspace{-0.16cm}
$\tau_{e}$ &exhaust gas time constant, [$s$]\\ %\vspace{-0.16cm}
$\rho$ & tunable DSMC gain, [-]\\
%\vspace{-0.16cm}
$T$ &sampling time, [$s$]\\ %\vspace{-0.16cm}
%$\tau_{e}$ &exhaust gas time constant, [$s$]\\ %\vspace{-0.16cm}
$\rho_{\beta}$ &adaptation gain for multiplicative model uncertainty, [$-$]\\ %\vspace{-0.16cm}
$\kappa$ &adaptation gain for additive model uncertainty, [$-$]\\ %\vspace{-0.16cm}
\end{tabular*}
\end{nomenclature}
\vspace{-1.0cm}
%%%%%%%%%%%%%%%%%%%%%%%%%%%%%%%%%%%%%%%%%%%%%%%%%%%%%%%%%%%%%%%%%%%%%%
\section{INTRODUCTION} \label{Sec:Intro}
Minimizing the gap between the designed and the implemented controllers is a major challenge in the design cycle of model-based industrial controllers. This gap is created due to (i) digital implementation of controller software that introduces sampling and quantization imprecisions at controller input/output (I/O) via analog to digital conversion (ADC), and (ii) uncertainty in modeling the plant dynamics. Model uncertainties are due to physical changes in the plant (e.g., aging) or fluctuations in the environment in which the system operates (e.g., temperature and humidity variations). %directly propagates through the controller by implementing the hardware. 
These sources of uncertainties should be addressed in the early stages of the controller design cycle, otherwise they substantially affect controller performance~\cite{NASA_Dabney,shahbakhti2012early}, and consequently increase the time and cost of verification and validation (V\&V) iterations.

\begin{figure}[h!]
\begin{center}
\includegraphics[angle=0,width= \columnwidth]{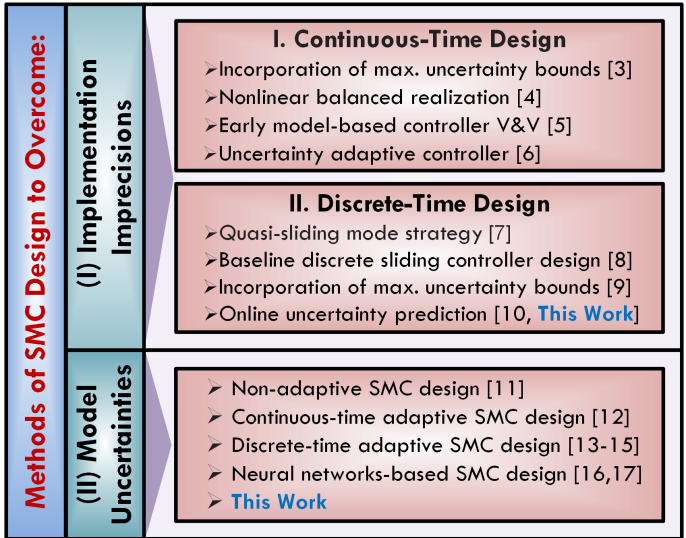} \vspace{-0.75cm}
\caption{\label{fig:Survay_Sampling_DSCC}Background of previous studies\cite{KyleACC,Amini_SAE,Amini_DSC,AminiSAE2016,Hansen_DSCC,Kyle_CDC,Kyle_DSCC,Amini_ACC2016,Misawa_DSC,Slotine,Acary,Chan_Automatica,Selina_PhD,munoz2000adaptive,fang1999use} on controller design
against ADC and modeling uncertainties.} \vspace{-1cm}
\end{center}
\end{figure}

Previous studies in the literature (shown in Fig.~\ref{fig:Survay_Sampling_DSCC}) have demonstrated the capabilities of sliding mode controller (SMC) in handling implementation and model uncertainties. The SMC design with robustness to implementation imprecisions can be divided into two major groups: (i) continuous-time designs, and (ii) discrete-time designs. For the first group, results in \cite{KyleACC} showed that by incorporating the maximum ADC uncertainty bounds on state equations, the robustness of the implemented controller improves against implementation uncertainties. The other effective approaches to enhance the controller robustness characteristics are nonlinear balanced realization~\cite{Amini_SAE}, and examining the controller structure to find states/variables with the largest numerical noise and then, removing those noise sources from the controller structure ~\cite{Amini_DSC}. The last approach is to incorporate an uncertainty adaptation mechanism into the SMC formulation by which the control actions can be modified to compensate the impact of ADC uncertainties~\cite{AminiSAE2016}.\vspace{-0.5cm}

Our previous works showed that a discrete sliding mode controller (DSMC), which explicitly includes the sampling time, shows better performance under sampling imprecisions compared to a continuous-time SMC that is digitally implemented~\cite{Hansen_DSCC}. Quasi and baseline (no incorporated ADC imprecisions) DSMCs under implementation imprecisions were introduced in~\cite{Hansen_DSCC} and~\cite{Kyle_CDC}, respectively. The robustness of the DSMC design was improved by incorporating the maximum ADC uncertainty bounds on control inputs in~\cite{Kyle_DSCC}. The results showed that the proposed DSMC with maximum ADC uncertainty bounds can result in better tracking performance compared to a baseline SMC/DSMC. The incorporation of maximum ADC uncertainty bounds requires an ideal symbolic model of the controller, which prevents establishment of an online framework to calculate uncertainty bounds. Moreover, the results in~\cite{KyleACC} and~\cite{Kyle_DSCC} showed that the calculated maximum ADC uncertainty bounds, which are based on a worst case scenario, leads to a conservative controller design and also large control actions. To this end, in our previous work~\cite{Amini_ACC2016}, we developed an online technique to predict and propagate sampling and quantization imprecisions on control signals. The results showed that not only is the conservative controller design avoided when using the new ADC uncertainty prediction technique, but also that the controller robustness improves significantly compared to a baseline DSMC which includes no uncertainty knowledge inside its structure.      

Handling the modeling uncertainty in sliding mode controller design has been the subject of several previous works in the literature. The model uncertainty/mismatch has been addressed using a non-adaptive framework~\cite{Misawa_DSC}, continuous-time adaptive form~\cite{Slotine}, discrete-time adaptive form~\cite{Acary,Chan_Automatica,Selina_PhD}, and neural network-based SMC designs~\cite{munoz2000adaptive,fang1999use}. The adaptive DSMC formulation from~\cite{Selina_PhD} presents a generic solution for compensating multiplicative type of model uncertainty for a general class of affine nonlinear systems. However, the proposed adaptive DSMC was not investigated under sampling and quantization imprecisions. Besides, regardless of the state equations of the physical model, it was assumed that an unknown multiplicative term can represent the existing uncertainty in the model.

The contribution of this paper is twofold. First, an online technique is presented to predict and propagate the ADC uncertainties on a measured signal based on a modified DSMC formulation. Second, the modified DSMC with predicted implementation imprecisions is extended to handle the additive type of model uncertainty using a discrete Lyapunov stability analysis. The proposed adaptive DMSC with predicted implementation imprecisions provides an integrated framework to improve the robustness of a baseline DSMC against both model and hardware (ADC) uncertainties. The proposed adaptive DMSC is verified on a real ECU.\vspace{-0.5cm}
\section{Incorporation of ADC Uncertainties into DSMC} \vspace{0.1cm}
\label{sec:UncertaintyPrediction}
\subsection{Uncertainty Prediction on Measured Signals} \label{sec:ADCUncertaintyPrediction}
\begin{figure}[h!]
\begin{center}
\includegraphics[angle=0,width= \columnwidth]{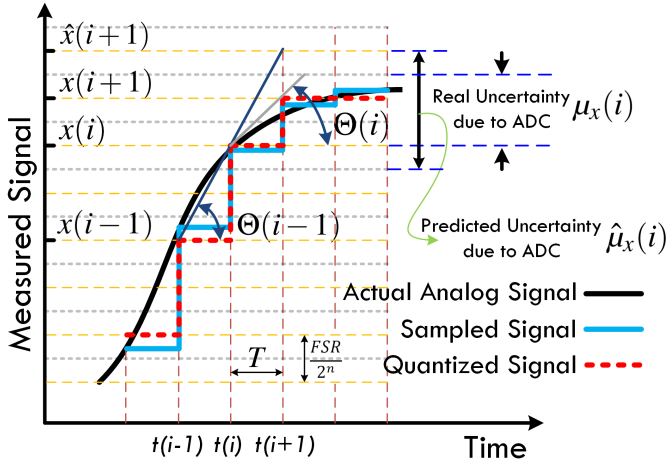} \vspace{-1cm}
\caption{\label{fig:ADC_Imprecisions}Uncertainty on a measured analog signal due to sampling and quantization.} \vspace{-0.75cm}
\end{center}
\end{figure}

Fig.~\ref{fig:ADC_Imprecisions} shows the ADC imprecisions introduced by sampling and quantization on a measured analog signal. According to the proposed uncertainty prediction method in~\cite{Amini_ACC2016}, the ADC uncertainty on a measured signal can be estimated using the slope of the signal and the introduced round-off error due to quantization. As can be seen in Fig.~\ref{fig:ADC_Imprecisions}, the sampling uncertainty $\left(\mu_{x_s}(i)\right)$ on a measured signal at each time step ($i$) can be obtained with respect to the signal slope ($\Theta$) as follows: \vspace{-0.4cm}
\begin{gather}\label{eq:StageIII_1}
\mu_{x_s}(i)={\Theta}(i)\times T 
\end{gather}
where $T$ is sampling time. Ideally, ${\Theta}(i)$ is calculated with respect to the value of the sampled signal at the current and next time steps. However, in practice, the sampled data are digitized at a rate equal to the ADC quantization level. Therefore, $\Theta$ is calculated with respect to the sampled and quantized signal:
\vspace{-0.4cm}
\begin{gather}\label{eq:StageIII_2}
{\Theta}(i)=\frac{x(i+1)-x(i)}{T}
\end{gather}
\vspace{-0.65cm}

In Eq.~(\ref{eq:StageIII_2}), the value of $x(i+1)$ is unknown at $i^{th}$ time step, thus its value is predicted $\left(\hat{x}(i+1)\right)$ to obtain the uncertainty due to ADC. Here, the slope from the previous time instant ($i-1$) is calculated and it is assumed that the slope at $i^{th}$ time step (${\hat{\Theta}}(i)$) is equal to the previous time instant $\left({\Theta}(i-1)\right)$: \vspace{-0.3cm}
%\vspace{-1cm}
\begin{gather}\label{eq:StageIII_3}
\hat{\Theta}(i)=\Theta(i-1) \Rightarrow \hat{x}(i+1)=T \times {\hat{\Theta}}(i) + x(i)
\end{gather}
%\vspace{-0.1cm}
%which leads to:\vspace{-0.4cm}
%\vspace{-1cm}
%\begin{gather}\label{eq:StageIII_4}
%\hat{x}(i+1)=2x(i)-x(i-1)
%\end{gather}
Moreover, the uncertainty introduced by quantization $\left(\mu_{x_q}(i)\right)$ can be calculated at each time step by considering a round-off error~\cite{Amini_DSC}:
\vspace{-0.5cm}
\begin{gather}\label{eq:StageIII_6}
\mu_{x_q}(i)=\frac{1}{2} \frac{FSR}{2^n}
\end{gather}
where $FSR$ is the full scale range of the measured signal and $n$ is the ADC's number of bits and represents the ADC resolution. Overall, the uncertainty on measured signals because of sampling and quantization at each time step can be expressed as follows using Eq.~(\ref{eq:StageIII_1}) and~(\ref{eq:StageIII_6}):\vspace{-0.5cm}
\begin{gather}\label{eq:StageIII_7}
\mu_{x}(i)=\hat{\mu}_{x_s}(i)+\mu_{x_q}(i)=x(i)-x(i-1)+\frac{1}{2} \frac{FSR}{2^n}
\end{gather}
where $\hat{\mu}_{x_s}$ is the predicted uncertainty due to sampling with respect to the value of the predicted sampled signal at $i+1^{th}$ time step ($\hat{x}_{i+1}$).
%The accuracy of~(\ref{eq:StageIII_7}) is sensitive to the magnitude of the changes in the signal's slope. Moreover, the accuracy of the uncertainty estimation depends on the ADC quantization level. Lower ADC resolution leads to poor uncertainty estimation. During the time interval between two time steps, output of a controller is constant, therefore the predicted value for ADC imprecisions at each time step can be assumed to be constant over the next time interval. In other words, 
Eq.~(\ref{eq:StageIII_7}) determines ADC uncertainty prediction at each time step until the next time step, since during the time interval between two time steps, output of a controller is constant. \vspace{-0.5cm}

\subsection{Uncertainty Propagation via DSMC Design} \label{sec:ControlUncertaintyPrediction}
A general class of single-input single-output (SISO) affine nonlinear systems can be expressed as follows:
\vspace{-0.37cm}
\begin{gather}\label{eq:GeneralNonlinear}
\dot{x}(t)=f(x(t))+g(u(t))
\end{gather}
where $x\subset{X\in{\mathbb{R}^{w}}}$, and $u\subset{U\in{\mathbb{R}^{v}}}$ are the state, and the input variables, respectively. $g$ is a non-singular square matrix and $f(x(t))$ represents the dynamics of the plant and does not depend on the inputs. For the system described in Eq.~(\ref{eq:GeneralNonlinear}), a SISO discrete sliding mode controller (DSMC) is designed to drive the state of the system ($x$) to its desired value ($x_d$). DSMC converts a complicated dynamical system into a first-order system through a sliding surface
transformation. Here, a DSMC design approach is adopted from~\cite{Kyle_CDC,Amini_ACC2016}. The sliding surface variable ($s$) is defined as the difference between desired ($x_d$) and measured signal ($x$) at each time step as follows:
\vspace{-0.4cm}
\begin{gather}\label{eq:StageIII_9}
%\begin{cases}
s(i)=x(i)-x_d(i)
%S(i+1)=X(i+1)-X_d(i+1)
%\end{cases}
\end{gather}
%\vspace{-0.75cm}
%
The control input $u(i)$ is obtained according to the following sliding reaching law:
\vspace{-0.5cm}
\begin{gather}\label{eq:StageIII_10}
|s(i+1)|\leq \rho |s(i)|
\end{gather}
where $0<\rho<1$ is the tunable DSMC gain. The first order Euler approximation, $\text{$\dot{x}\approx T^{-1}\times({x(i+1)-x(i)}$})$, is used to discretize the nonlinear system in Eq.~(\ref{eq:GeneralNonlinear}):
\vspace{-0.35cm}
\begin{gather}\label{eq:StageIII_16}
%\dot{X}=f(X(t))+gU(t) \xrightarrow{\text{$\dot{x}=\frac{x(i+1)-x(i)}{T}$}}\\ \nonumber
x(i+1)=x(i)+fT+gu(i)T
\end{gather}
%\vspace{-0.75cm}
%

\begin{figure}[h!]
\begin{center}
\includegraphics[width=\columnwidth]{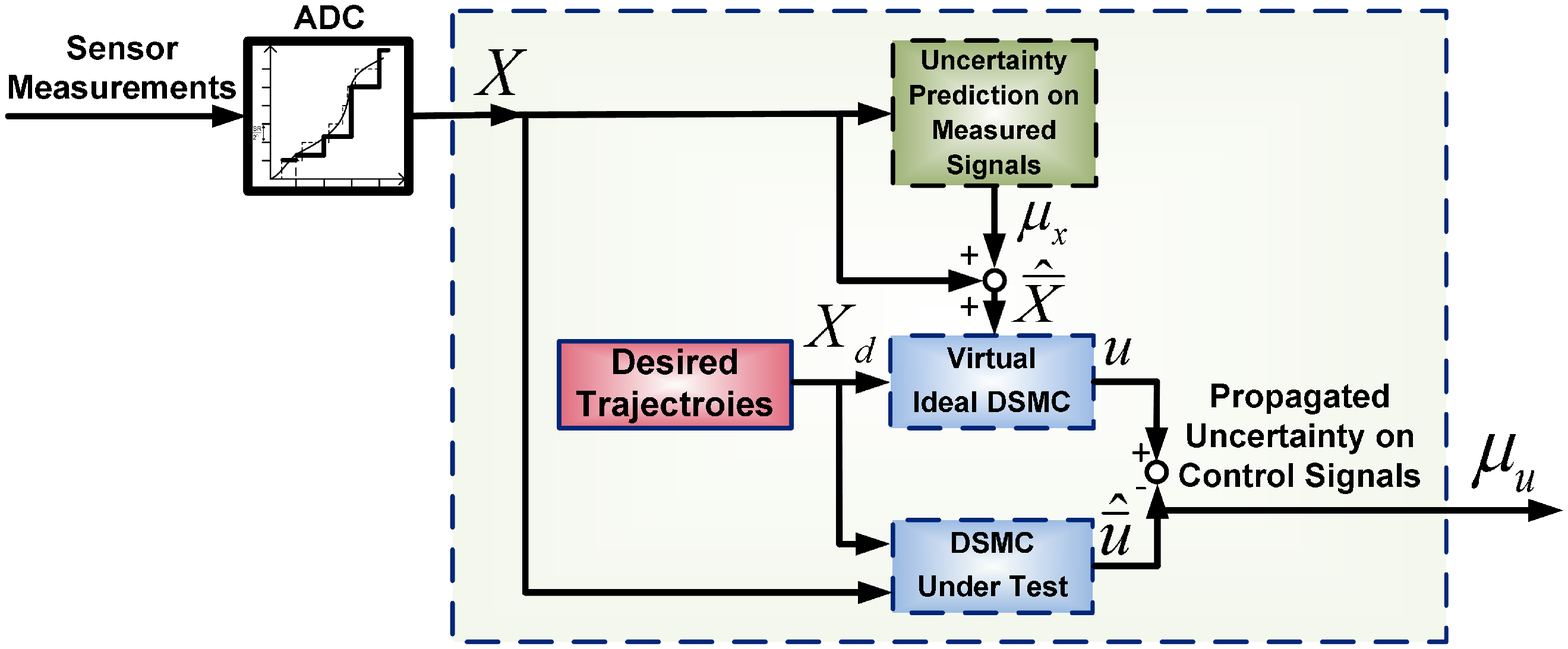} \vspace{-0.6cm}
\caption{\label{fig:dsmc_with_PU}Schematic of the DSMC with online estimations of uncertainties on measured and control signals.} \vspace{-0.8cm}
\end{center}
\end{figure}

The baseline DMSC can be formulated for the discrete nonlinear system in Eq.~(\ref{eq:StageIII_16}) with the sliding surface function and reaching law described in Eq.~(\ref{eq:StageIII_9}) and~(\ref{eq:StageIII_10}) as follows:
\vspace{-0.35cm}
\begin{gather}\label{eq:StageIII_17}
u(i)=\frac{1}{gT}[(\rho-1)x(i)-\rho x_d(i)-fT+x_d(i+1)]
\end{gather}
%\vspace{-0.35cm}
%
According to the proposed approach in~\cite{Amini_ACC2016,Kyle_DSCC}, the structure of the baseline DSMC in Eq.~(\ref{eq:StageIII_17}) is
modified against ADC uncertainties by incorporating the propagated implementation imprecisions on control signals into the DSMC formulation. To this end, the ADC uncertainties on the measured signals should be first estimated, and then be propagated on control signals. Unlike linear systems, in nonlinear systems it is not possible to conduct the noise propagation on control signals analytically~\cite{Amini_ACC2016}. Fig.~\ref{fig:dsmc_with_PU} shows our proposed approach to estimate the propagated uncertainties on control signals for the nonlinear system described in Eq.~(\ref{eq:StageIII_16}). The uncertainty prediction mechanism contains two auxiliary baseline DSMCs, including a virtual ideal DSMC and the DSMC under test. Desired trajectories for both auxiliary DSMCs are the same. The feedback signals to the DSMC under test are vectors of the measured signals after ADC ($X$), while the inputs to the virtual ideal DSMC are the estimated vector of the measured signals ($\hat{\bar{X}}$) before ADC. $\hat{\bar{X}}$ is estimated according to the real values of the measured signals after ADC and the predicted ADC uncertainties on measured signals ($\mu_{_X}=diag[\mu_{x_1},\mu_{x_2},...,\mu_{x_w}]$) from Eq.~(\ref{eq:StageIII_7}):\vspace{-0.35cm}
\begin{gather}\label{eq:StageIII_17_2}
\hat{\bar{X}}=\mu_{_X}+X
\end{gather}
\vspace{-0.8cm}
%
%where $\mu_{_X}=diag[\mu_{x_1},\mu_{x_2},...,\mu_{x_w}]$. %\vspace{-0.45cm}
%

Since the only difference between the DSMC under test and the virtual ideal DSMC is the ADC uncertainties on measured signals, the propagated uncertainty vector on control signals is the difference between control signals of these two controllers as follows:
\vspace{-0.5cm}
\begin{gather}\label{eq:StageIII_18}
\mu_{_U}=diag[\hat{\bar{u}}_1-u_1,...,\hat{\bar{u}}_{v}-u_{v}]
\end{gather}
%\vspace{-0.75cm}
where $\hat{\bar{u}}_{1,...,{v}}$ are the estimated control inputs from the ideal DSMC ($\bar{u}$) according to the predicted ADC uncertainties on measured signals. \vspace{-0.35cm}
\subsection{Inclusion of ADC Uncertainties into the DSMC}
Structure of the baseline DSMC (Eq.~(\ref{eq:StageIII_17})) is modified by inclusion of estimated uncertainties on control signals into the DSMC formulation:
\vspace{-0.4cm}
\begin{gather}\label{eq:StageIII_19}
u^{mod}(i)=u(i)+\delta(i)
\end{gather}
where $u^{mod}$ is the modified controller input and $u$ is calculated with respect to Eq.~(\ref{eq:StageIII_17}). $\delta(i)$ is included to compensate for hardware imprecisions. $\delta(i)$ is obtained according to the predicted uncertainties on control signals as follows:
\vspace{-0.4cm}
\begin{gather}\label{eq:StageIII_20}
\delta(i)=-\mu_{u}(i)s(i)
\end{gather}
%\vspace{-0.75cm}
By inclusion of propagated ADC uncertainties on control signals ($\mu_{u}$), the baseline DSMC for a nonlinear system in Eq.~(\ref{eq:StageIII_16}) is modified as follows:
\vspace{-0.4cm}
\begin{gather}\label{eq:StageIII_22}
u^{mod}(i)=
\frac{1}{gT}[(\rho-1)x(i)-\rho x_d(i)-fT+x_d(i+1)]-\mu_{u}(i)s(i)
\end{gather} \vspace{-1.5cm}
\section{Model Uncertainty}
The uncertainties in the modeled dynamics ($f$) can be expressed using either multiplicative ($\beta$) or additive ($\alpha$) terms in the presence of hardware (ADC) imprecisions ($\delta=-\mu_u\times s$) using the nominal discrete nonlinear model in Eq.~(\ref{eq:StageIII_16}): \vspace{-0.4cm}
\begin{subequations} \label{eq:Model_un1}
\begin{align}
x(i+1)=x(i)+T(\beta\times f)+Tg(u(i)+\delta(i))
\end{align} \vspace{-1.5cm}
\begin{align}
x(i+1)=x(i)+T(f+\alpha)+Tg(u(i)+\delta(i))
\end{align} %\vspace{-0.5cm}
\end{subequations}
The selection of the model uncertainty type depends on the physics and state equations. In the following sections, the adaptation mechanism for the multiplicative type of model uncertainty (Eq.~(\ref{eq:Model_un1}a)) is first adopted from literature, and then it is combined with the ADC uncertainty prediction presented in the previous section to design an integrated adaptive DSMC that is robust to the multiplicative type of model uncertainty and implementation imprecisions. Next, the adaptation law for overcoming the additive type of model uncertainty (Eq.~(\ref{eq:Model_un1}b)) will be derived using Lyapunov stability theorem. \vspace{-0.45cm}
\subsection{Multiplicative Type of Model Uncertainty} \label{multiplicative_adaptive_section}
Here, an adaptation mechanism is adopted from~\cite{Selina_ACC2016} to address the multiplicative type of model uncertainty (Eq.~(\ref{eq:Model_un1}a)) by which the model with error is driven to the nominal model. The following adaptation law is used to update the unknown constant multiplicative term ($\beta$) in the model. The objective is converging $\beta$ to its nominal value ``1" from an unknown initial condition:  \vspace{-0.5cm}
\begin{gather}\label{eq:Model_un_selina1}
\hat{\beta}(i+1)=\hat{\beta}(i)+\frac{fs(i)T}{{\rho_{_\beta}}}
\end{gather}
where $\hat{\beta}$ is the estimation of the unknown parameter ($\beta$), $f$ is the nominal model, $s$ is the sliding function, and $\rho{_\beta}$ is the adaptation gain for the multiplicative type of model uncertainty. Eq.~(\ref{eq:StageIII_22}) with incorporated implementation imprecisions can be rewritten to include the multiplicative uncertainty term in the model as follows: \vspace{-0.6cm}
\begin{gather}\label{eq:Model_un_selina2}
u^{mod}_{adaptive,~multiplicative}(i)= \\ 
\frac{1}{gT}[(\rho-1)x(i)-\rho {x}{_d}(i)-(\hat{\beta} \times f)T+{x}{_d}(i+1)]-\mu_u(i)s(i) \nonumber
\end{gather}
By solving Eq.~(\ref{eq:Model_un_selina1}) and ({\ref{eq:Model_un_selina2}}) simultaneously, the adaptive DSMC with incorporated implementation imprecisions is able to (i) enhance controller robustness characteristics against ADC uncertainties compared to the baseline DSMC (Eq.~(\ref{eq:StageIII_17})), and (ii) remove the error in the modeled dynamics by driving the model with unknown multiplicative uncertainty to its nominal value. \vspace{-0.35cm}  
\subsection{Additive Type of Model Uncertainty}
\begin{figure}[h!]
\begin{center}
\includegraphics[width=\columnwidth]{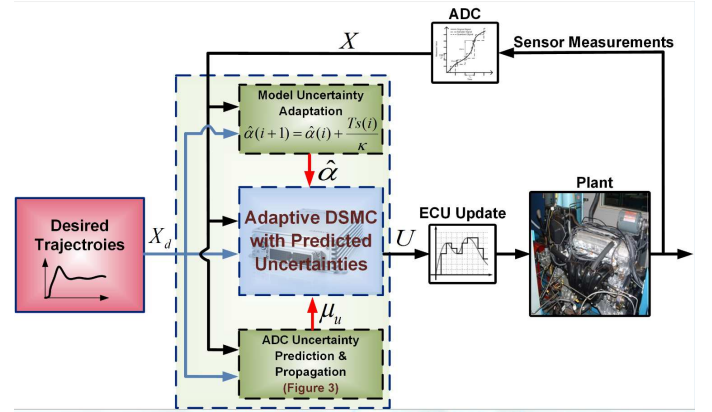} \vspace{-0.6cm}
\caption{\label{fig:dsmc_with_PU_full} Schematic of the adaptive DSMC with online estimations of uncertainties on measured and control signals.} \vspace{-1.5cm}
\end{center}
\end{figure}

In this section, the adaptive DSMC theory is developed for addressing the additive type of model uncertainty using Lyapunov stability theory. By expanding the sliding reaching law in Eq.~(\ref{eq:StageIII_10}) for the discrete nonlinear system in (Eq.~(\ref{eq:Model_un1}b)) with additive uncertainty term ($\alpha$), the following relationship is concluded: \vspace{-0.4cm}
\begin{gather}\label{eq:Model_un_mamini1}
s(i+1)=x(i)+(f+\alpha)T+g u(i)-{{x}}{_d}(i+1)
\end{gather}
where $\alpha$ is unknown and constant. Control input ($u$) is calculated according to Eq.~(\ref{eq:StageIII_17}) by inclusion of the uncertainty term:\vspace{-0.4cm}
\begin{gather}\label{eq:Model_un_mamini2}
u(i)=
\frac{1}{g T}[(\rho-1)x(i)-\rho {x}{_d}(i)-(f+\hat{\alpha}(i))T+{x}{_d}(i+1)]
\end{gather}%\vspace{-1cm}
where $\hat{\alpha}$ is the estimation of the unknown parameter. By substituting Eq.~(\ref{eq:Model_un_mamini2}) into Eq.~(\ref{eq:Model_un_mamini1}), the sliding surface dynamics can be simplified as follows: \vspace{-0.4cm}
\begin{gather}\label{eq:Model_un_mamini3}
s(i+1)=\rho s(i)+T(\alpha-\hat{\alpha}(i))
\end{gather}%\vspace{-1cm}
A new term ($\tilde{\alpha}$) is defined to represent the difference between the unknown and estimated additive uncertainty term ($\tilde{\alpha}(i)=\alpha-\hat{\alpha}(i)$). Upon substitution of $\tilde{\alpha}$ into Eq.~(\ref{eq:Model_un_mamini3}), the sliding surface dynamics in terms of the DSMC gain ($\rho$) and the error in the estimation of the unknown parameter ($\tilde{\alpha}$) can be expressed as: \vspace{-0.4cm}
\begin{gather}\label{eq:Model_un_mamini4}
s(i+1)=\rho s(i)+T\tilde{\alpha}(i)
\end{gather}
A Lyapunov-based analysis is employed to first determine the stability of the closed-loop system, and second, derive the adaptation law. To this end, the following Lyapunov function candidate is proposed: \vspace{-0.4cm}
\begin{gather}\label{eq:Model_un_mamini5}
V(i)=\frac{1}{2} {s}^2(i)+\frac{1}{2}\kappa{\tilde{\alpha}}^2(i)
\end{gather}
where $\kappa>0$ is a tunable parameter chosen for the numerical sensitivity of the unknown parameter estimation (the adaptation gain for the additive type of model uncertainty). The proposed Lyapunov function is positive definite and quadratic with respect to the sliding variable ($s(i)$) and the additive unknown parameter estimation error ($\tilde{\alpha}$). Similar to continuous-time systems, in which the negative definite condition is required for the derivative of the Lyapunov function to guarantee the asymptotic stability, in the discrete time domain, investigation of a Lyapunov difference equation is required for the Lyapunov stability analysis. The Lyapunov difference function can be calculated using an implicit approach in which the value of the Lyapunov function at the subsequent time step is first obtained using a two-variable Taylor series expansion:  \vspace{-0.4cm}
\begin{gather}\label{eq:Model_un_mamini6}
V(i+1)=V(i)+ \frac{\partial V(i)}{\partial s(i)}\Delta s(i)+ \frac{\partial V(i)}{\partial \tilde{\alpha}(i)}\Delta \tilde{\alpha}(i)+\\ \nonumber
\frac{1}{2} \frac{\partial^2 V(i)}{\partial {s}^2(i)}\Delta {s}^2(i)+
\frac{1}{2} \frac{\partial^2 V(i)}{\partial {\tilde{\alpha}}^2(i)}\Delta {\tilde{\alpha}}^2(i)+ \\ \nonumber
\frac{\partial^2 V(i)}{\partial {s}(i)\partial \tilde{\alpha}(i)}\Delta {s}(i)\times \Delta \tilde{\alpha}(i)+...
\end{gather}\vspace{-0.25cm}
where,\vspace{-0.4cm}
\begin{gather}\label{eq:Model_un_mamini7}
\Delta s(i)\equiv s(i+1)-s(i)\\ \nonumber
\Delta \tilde{\alpha}(i) \equiv \tilde{\alpha}(i+1)-\tilde{\alpha}(i)
\end{gather}
%
%and the partial derivatives are as follows: \vspace{-0.25cm}
%
%\begin{gather}\label{eq:Model_un_mamini8}
%\frac{\partial V(i)}{\partial s(i)}=s(i),~~\frac{\partial V(i)}{\partial \tilde{\alpha}(i)}=\kappa\tilde{\alpha}(i) \\ \nonumber
%\frac{\partial^2 V(i)}{\partial {s}^2(i)}=1,~~\frac{\partial^2 V(i)}{\partial {\tilde{\alpha}}^2(i)}=\kappa,~~\frac{\partial^2 V(i)}{\partial {s}(i)\partial \tilde{\alpha}(i)}=0
%\end{gather}%\vspace{-1cm}
%
Next, the Lyapunov difference function ($\Delta V(i)=V(i+1)-V(i)$) is calculated by substituting the values of the partial derivatives into Eq.~(\ref{eq:Model_un_mamini6}): \vspace{-0.4cm}
\begin{gather}\label{eq:Model_un_mamini9}
\Delta V(i) 
=s(i)\Delta s(i)+ \kappa\tilde{\alpha}(i)\Delta \tilde{\alpha}(i)+
\frac{1}{2} \Delta {s}^2(i)+ \frac{1}{2} \kappa\Delta {\tilde{\alpha}}^2(i)+...
\end{gather}
As can be observed from Eq.~(\ref{eq:Model_un_mamini9}), all higher order ($>2$) derivatives are zero. After substitution Eq.~(\ref{eq:Model_un_mamini7}) in Eq.~(\ref{eq:Model_un_mamini9}), the Lyapunov difference function is simplified as follows: \vspace{-0.4cm} 
%
%\begin{gather}\label{eq:Model_un_mamini10}
%\Delta V_k(i)=s_k(i)(s_k(i+1)-s_k(i))+ \kappa_k\tilde{\alpha}_k(i)(\tilde{\alpha}_k(i+1)-\tilde{\alpha}_k(i))+ \nonumber \\
%\frac{1}{2} \Delta {s_k}^2(i)+ \frac{1}{2} \kappa_k\Delta {\tilde{\alpha}}^2_k(i)+...
%\end{gather}\vspace{-1cm}
%
%\begin{gather}\label{eq:Model_un_mamini11}
%\Delta V_k(i)=s_k(i)(\rho_k s_k(i)+T\tilde{\alpha}_k(i)-s_k(i))+ \nonumber \\ \kappa_k\tilde{\alpha}_k(i)(\tilde{\alpha}_k(i+1)-\tilde{\alpha}_k(i))+ \nonumber \\
%\frac{1}{2} \Delta {s_k}^2(i)+ \frac{1}{2} \kappa_k\Delta {\tilde{\alpha}}^2_k(i)+...
%\end{gather}\vspace{-1cm}
%
\begin{gather}\label{eq:Model_un_mamini12}
\Delta V(i)=(\rho-1){s}^2(i)+T\tilde{\alpha}(i)s(i)+ \nonumber \\ \kappa\tilde{\alpha}(i)(\tilde{\alpha}(i+1)-\tilde{\alpha}(i))+ 
\frac{1}{2} \Delta {s}^2(i)+ \frac{1}{2} \kappa\Delta {\tilde{\alpha}}^2(i)
\end{gather}
which yields: \vspace{-0.4cm}
\begin{gather}\label{eq:Model_un_mamini13}
\Delta V(i)=-(1-\rho){s}^2(i)+\kappa \tilde{\alpha}(i) \left (\tilde{\alpha}(i+1)-\tilde{\alpha}(i)+\frac{Ts(i)}{\kappa} \right) \nonumber  \\ +O({\Delta {s}^2(i),\Delta\tilde{\alpha}}^2(i))
\end{gather}
The following adaptation law is chosen to update the error in unknown parameter estimation $\tilde{\alpha}$: \vspace{-0.4cm}
\begin{gather}\label{eq:Model_un_mamini14}
\tilde{\alpha}(i+1)= \tilde{\alpha}(i)-\frac{Ts(i)}{\kappa} 
\end{gather}
Assuming that for small enough sampling periods $O({\Delta {s}^2(i),\Delta\tilde{\alpha}}^2(i))\approx 0$~\cite{Selina_PhD}, upon incorporating the adaptation law from Eq.~(\ref{eq:Model_un_mamini14}) in the Lyapunov difference function, we have: \vspace{-0.5cm}
\begin{gather}\label{eq:Model_un_mamini15}
\Delta V(i)=-(1-\rho){s}^2(i)
\end{gather}
It can be concluded from Eq.~(\ref{eq:Model_un_mamini15}) that there exists a region around $s(i)=0$ and $\tilde{\alpha}=0$ where the Lyapunov difference function is negative semi-definite. This means that the sliding variable (the tracking error, $s$) converges to zero, and the error in estimating the unknown additive parameter ($\tilde{\alpha}$) is at least bounded~\cite{Selina_ACC2016}. In the next step, the adaptive baseline DSMC (Eq.~(\ref{eq:Model_un_mamini2})), in which the unknown additive error in the model ($\hat{\alpha}$) is driven to zero by solving the adaptation law (Eq.~(\ref{eq:Model_un_mamini14})), is modified against implementation impressions by incorporating the knowledge of ADC imprecisions on control signals according to Eq.~(\ref{eq:StageIII_19}): \vspace{-0.6cm}
%
%\begin{gather}\label{eq:Model_un_mamini16}
%\hat{\alpha}_k(i+1)=\hat{\alpha}_k(i)+\frac{Ts_k(i)}{\kappa_k}
%\end{gather}\vspace{-1cm}
%
\begin{gather}
\label{eq:Model_un_mamini17}
u^{mod}_{adaptive,~additive}(i)=\\ 
\frac{1}{gT}[(\rho-1)x(i)-\rho {x}{_d}(i)-(\hat{\alpha} + f)T+{x}{_d}(i+1)]-\mu_u(i)s(i) \nonumber 
\end{gather}
Fig.~(\ref{fig:dsmc_with_PU_full}) shows the schematic of the proposed integrated adaptive DSMC. In the next section, the new adaptive DSMC is tested on an engine control problem under implementation and additive and multiplicative types of model uncertainties. \vspace{-0.7cm}
\section{Case Study: Engine Control during Cold Start}\label{subsec:Engine_Control}
The effectiveness of adaptive control stategies for engine and powertrain controls has been demonstrated in the literature~\cite{JingSun_1,JingSun_2}. Here, a physics-based spark ignition (SI) engine model during cold start from \cite{Shaw} is used to investigate the proposed adaptive DSMC performance. The model is parameterized for a 2.4-liter, 4-cylinder, DOHC 16-valve Toyota 2AZ-FE engine and a three-way catalyst exhaust aftertreatment system. The engine rated power is 117kW $@$ 5600~RPM and it has a rated torque of 220 Nm $@$ 4000~RPM. The experimental validation of different components of the engine model is found in~\cite{Sanketi,Salehi_ACC}. The model has three inputs: (i) $\dot{m}_{ai}$, air mass flow rate into the intake manifold, (ii) $\dot{m}_{fc}$, commanded fuel mass flow rate, and (iii) $\Delta$, spark timing after top dead center (ATDC). The discretized model includes four difference equations:
\vspace{-0.40cm}
%\begin{gather}
%\label{eq:ModelStates_1}\dot{x}_1=\dot{m}_a= \dot{m}_{ai} - \dot{m}_{ao}(m_a,\omega_e) \\
%\label{eq:ModelStates_2}\dot{x}_2= \dot{\omega}_e= \frac{1}{J}\left[T_E(m_a,\omega_e)\right]\\
%\label{eq:ModelStates_3}\dot{x}_3= \ddot{m}_f= \frac{1}{\tau_f}\left[\dot{m}_{fc} - \dot{m}_{f}\right]\\
%\label{eq:ModelStates_4}\dot{x}_4= \dot{T}_{cat}= \frac{1}{mC_p}\left[\dot{Q}_{gen} + \dot{Q}_{in} - \dot{Q}_{out}\right]\\
%\label{eq:ModelStates_5}\dot{x}_5=\dot{T}_{exh}=\frac{1}{\tau_e}\left[SI(\Delta).AI(m_a,\omega_e,\dot{m}_f)-T_{exh}\right]
%\end{gather} \vspace{-0.5cm}
%
\begin{gather}
\label{eq:Engine_discretized}
T_{exh}(i+1)=(1-\frac{T}{\tau_e})T_{exh}(i)+\frac{T}{\tau_e}(7.5\Delta(i)+600)AFI(i)  \\ 
\dot{m}_f(i+1)=\dot{m}_f(i)+\frac{T}{\tau_f}[\dot{m}_{fc}(i)-\dot{m}_f(i)] \nonumber \\ 
\omega_e(i+1)=\omega_e(i)+\frac{T}{J}T_E(i) \nonumber \\ 
m_a(i+1)=m_a(i)+[\dot{m}_{ai}(i)-\dot{m}_{ao}(i)]T \nonumber
\end{gather} 
The four states are the exhaust temperature~($T_{exh}$), fuel mass flow rate into the cylinders (${\dot{m}_f}$), the engine speed~(${\omega_e}$), and the mass of air inside the intake manifold ($m_{a}$). Details of the functions and constants in the engine model are found in ~\cite{Sanketi}. Here, the control objective is to keep exhaust temperature, engine speed, and air to fuel ratio (AFR) at their desired values. To this end, a set of four SISO DSMCs is designed. For better controller evaluation, desired engine trajectories are defined in their worst and non-smooth shapes. %\vspace{-0.1cm}

Each SISO DSMC is formulated according to a sliding surface. The first controller is the $T_{exh}$ controller for which the first sliding surface is defined to be the error in tracking the desired exhaust temperature, which is controlled by spark timing. The second sliding variable is defined to be the error in tracking the desired fuel mass flow ($\dot{m}_f$). Desired fuel mass flow ($\dot{m}_{f,d}$) is calculated from the desired AFR ($AFR_d=\frac{\dot{m}_{ao,d}}{\dot{m}_{f,d}}$).%, since AFR is a function of both air and fuel mass flow rates and air mass flow rate controls engine speed, the desired fuel mass flow rate is  and is considered as the second sliding surface. %Moreover, fuel mass flow rate is among the measured states of the model/engine which means the effect of the ADC imprecisions can be investigated on the designed SMC based on the assigned desired trajectory for the fuel mass flow rate. On the other hand, since AFR is calculated according to measured fuel and air mass flow rates after ADC, the sampling and quantization influences on AFR  cannot be explored explicitly, since those imprecisions affect both air and fuel mass flow rates.
~The third desired trajectory is defined to be the desired engine speed. Since there is no explicit control input for tracking the desired engine speed, desired air mass is considered as the synthetic control input for the engine speed regulation. Introducing the desired air mass requires another SMC to track desired air mass which means four sliding surfaces are required: \vspace{-0.45cm}
\begin{gather}\label{eq:ColdStart_SISO_SMC_SS}
s_{1}=T_{exh}-{T_{exh,d}},~~s_{2}={\dot{m}_{f}}-{\dot{m}_{f,d}}\\ \nonumber
s_3=\omega_e-{\omega_{e,d}},~~ s_4=m_a-{m_{a,d}}
\end{gather} %\vspace{-0.5cm}

$T_{exh}$, $AFR$ ($\dot{m}_f$), $m_a$, and $\omega_e$ can be measured on an engine; thus the ADC affects these four measured states and corresponding sliding surfaces. In the following sections, to drive the states of the system towards the desired values, a set of adaptive SISO DSMCs %(Eq.~(\ref{eq:Model_un_mamini2}) and (\ref{eq:Model_un_mamini14})) 
with incorporated ADC uncertainties (Eq.~(\ref{eq:Model_un_mamini14},\ref{eq:Model_un_mamini17})) is designed and will be discussed. %. Next, the predicted ADC uncertainties are incorporated into the adaptive DSMC formulation .  %Discrete equations of the engine model are required in the next step to design the DSMC according to Eq.~(\ref{eq:StageIII_16}). Discrete state equations of the engine model are shown below:
%
%\begin{gather}\label{eq:Engine_Disc}
%\begin{split}
%T_{exh}(i+1)=T_{exh}(i)+\frac{T}{\tau_e}[(7.5\Delta(i)+600)AFI(i)-T_{exh}(i)] \nonumber \\ 
%\end{gather}\vspace{-0.6cm}
%\begin{gather}
%\dot{m}_f(i+1)=\dot{m}_f(i)+\frac{T}{\tau_f}[\dot{m}_fc(i)-\dot{m}_f(i)] \nonumber \\ 
%\end{gather}\vspace{-0.6cm}
%\begin{gather}
%\omega_e(i+1)=\omega_e(i)+\frac{T}{J}T_E(i) \nonumber \\ 
%\end{gather}\vspace{-0.6cm}
%\begin{gather}
%m_a(i+1)=m_a(i)+[\dot{m}_ai(i)-\dot{m}_ao(i)]T 
%\end{split}
%\end{gather} \vspace{-0.5cm}
%Since no sliding surface was defined for the catalytic converter, its state equation (Eq.~(\ref{eq:ModelStates_4})) is not discretized. 
%

\vspace{0.25cm}
$\bullet { \textbf{~~Exhaust~Gas~ Temperature~Controller:}}$ According to Eq.~(\ref{eq:Engine_discretized}), the dynamics of the exhaust temperature ($f_{{T_{exh}}}$) with additive uncertainty term ($\alpha_{T_{exh}}$) is: \vspace{-0.45cm}
\begin{gather}
\label{eq:Engine_discretized_Texh}
f_{{T_{exh}}}=\frac{1}{\tau_e}[600AFI-T_{exh}+\alpha_{T_{exh}}]
\end{gather}
%
%Although in the empirical model of the exhaust gas temperature 
$\alpha_{T_{exh}}$ represents the error in the empirical model of the exhaust gas temperature. The error in the modeled exhaust gas temperature dynamics is compensated using the following adaptation law with respect to Eq.~(\ref{eq:Model_un_mamini14}): \vspace{-0.4cm}
\begin{gather}
\label{eq:adaptive_Texh}
\hat{\alpha}_{T_{exh}}(i+1)=\hat{\alpha}_{T_{exh}}(i)+\frac{T(T_{exh}(i)-T_{{exh,d}}(i))}{\kappa_1}
\end{gather}
By incorporating the solution of Eq.~(\ref{eq:adaptive_Texh}) and the predicted implementation imprecisions ($\mu_{\Delta}$) into Eq.~(\ref{eq:Model_un_mamini17}), the modified adaptive DSMC for exhaust gas temperature becomes: \vspace{-0.4cm}
\begin{gather}\label{eq:Engine_DSMC_Final_1}
\Delta(i)=\frac{\tau_e}{7.5~AFI\,.\,T}[-\frac{T}{\tau_e}(600~AFI-T_{exh}(i)+\hat{\alpha}_{T_{exh}}(i))~~~~~~  \\
+(\rho_1-1)s_1(i)+T_{exh,d}(i+1)-T_{exh,d}(i)]-\mu_{\Delta}(i)s_1(i) \nonumber
\end{gather} %\vspace{-1cm}
If $\mu_{\Delta}=0$ the controller is an adaptive baseline DSMC. But when $\mu_{\Delta}$ is calculated according to the mechanism in Fig.~(\ref{fig:dsmc_with_PU}), the controller is an integrated adaptive DSMC with predicted ADC uncertainties.

\vspace{0.25cm}
$\bullet { \textbf{~~Fuel~Flow~Rate~Controller:}}$ As can be observed from Eq.~(\ref{eq:Engine_discretized}), the fuel flow dynamics ($f_{\dot{m}_f}$) with additive uncertainty term ($\alpha_{\dot{m}_f}$) is as follows: \vspace{-0.5cm}
\begin{gather}
\label{eq:Engine_discretized_mdotf}
f_{\dot{m}_f}=-\frac{1}{\tau_f}[\dot{m}_f(i)+\alpha_{\dot{m}_f}] 
\end{gather}
In practice, $\dot{m}_f$ is not measured directly and it is calculated according to the AFR sensor measurement and estimated air mass flow into the cylinder ($\dot{m}_{ao}$). $\alpha_{\dot{m}_f}$ represents the error in estimating fuel flow rate due to AFR measurement uncertainty and error in predicting the air mass flow. The adaptation law for $\alpha_{\dot{m}_f}$ is:  
\vspace{-0.4cm}
\begin{gather}
\label{eq:adaptive_mdotf}
\hat{\alpha}_{\dot{m}_f}(i+1)=\hat{\alpha}_{\dot{m}_f}(i)+\frac{T(\dot{m}_f(i)-\dot{m}_{f,d}(i))}{\kappa_2}
\end{gather}
where $\dot{m}_{f,d}$ is calculated according to desired AFR. The adaptive control law for fuel flow rate into the cylinder with additive model uncertainty and predicted ADC imprecisions ($\mu_{\dot{m}_{fc}}$) is: \vspace{-0.35cm}
\begin{gather}\label{eq:Engine_DSMC_Final_2}
\dot{m}_{fc}(i)=\frac{\tau_f}{T}[\frac{T}{\tau_f}(\dot{m}_f(i)+\hat{\alpha}_{\dot{m}_f}(i))+(\rho_2-1)s_2(i)~~~~~~~~  \\
+\dot{m}_{f,d}(i+1)-\dot{m}_{f,d}(i)]- \mu_{\dot{m}_{fc}}(i)s_2(i) \nonumber
\end{gather} \vspace{-1cm}

\vspace{0.35cm}
$\bullet { \textbf{~~Engine~Speed~Controller:}}$
$f_{\omega_e}$ for the engine with additive uncertainty ($\alpha_{\omega_e}$) is as follows: \vspace{-0.4cm}
\begin{gather}
\label{eq:Engine_discretized_we}
f_{\omega_e}=-\frac{1}{J}(T_{loss}+\alpha_{\omega_e})
\end{gather} 
where $T_{loss}=0.4\omega_e+100$. $T_{loss}$ represents the torque losses (e.g., due to friction) on the crankshaft. Thus, the additive uncertainty $\alpha_{\omega_e}$ represents the error in reading the torque map. $\alpha_{\omega_e}$ is driven to zero using the following adaptation law: \vspace{-0.5cm}
\begin{gather}
\label{eq:adaptive_we}
\hat{\alpha}_{\omega_e}(i+1)=\hat{\alpha}_{\omega_e}(i)+\frac{T(\omega_e(i)-\omega_{e,d}(i))}{\kappa_3}
\end{gather}
Finally, the control input ($m_{a,d}$) for engine speed regulation after incorporating the propagated ADC uncertainties ($\mu_{m_{a,d}}$) becomes: \vspace{-0.25cm}
\begin{gather}\label{eq:Engine_DSMC_Final_3}
m_{a,d}(i)=\frac{J}{30,000\,T}[\frac{T}{J}(100+0.4\omega_e(i)+\hat{\alpha}_{\omega_e}(i))\\\nonumber+(\rho_3-1)s_3(i)  
+\omega_{e,d}(i+1)-\omega_{e,d}(i)]-\mu_{m_{a,d}}(i)s_3(i) ~~~~
%\end{split}
\end{gather} \vspace{-1cm}

\vspace{0.35cm}
$\bullet { \textbf{~~Air~Mass~Flow~Controller:}}$
The intake air manifold dynamics is linked to the rotational dynamics through the calculated $m_{a,d}$. The calculated $m_{a,d}$ from Eq.~(\ref{eq:Engine_DSMC_Final_3}) is used as the desired trajectory to obtain $\dot{m}_{ai}$ as the control input of the intake air flow rate controller. The intake air manifold mass dynamics with additive model uncertainty ($\alpha_{m_a}$) is: \vspace{-0.4cm}
\begin{gather}
\label{eq:Engine_discretized_ma}
f_{m_a}=-(\dot{m}_{ao}(i)+\alpha_{m_a})
\end{gather}  \vspace{-0.25cm}
where~\cite{Shaw}: \vspace{-0.4cm}
\begin{gather}
\label{eq:Engine_discretized_ma2}
\dot{m}_{ao}=k_1\eta_{vol}m_a\omega_e \\
\label{eq:Engine_discretized_ma3}
\eta_{vol}=m_a^2(k_2{\omega_e^2}+k_3\omega_e+k_4)+m_a(k_5{\omega_e^2}+k_6\omega_e+k_7) \\ \nonumber 
+k_8{\omega_e^2}+k_9{\omega_e}+k_{10}% \nonumber
\end{gather}
$k_{1,2,...,10}$ are the empirical (curve fit) parameters of the volumetric efficiency ($\eta_{vol}$). As can be seen from Eq.~(\ref{eq:Engine_discretized_ma})-(\ref{eq:Engine_discretized_ma3}), the additive uncertainty term in the intake air manifold dynamics ($\alpha_{m_a}$) compensates for the uncertainty in $\dot{m}_{ao}$ that is read through $\eta_{vol}$ curve fit. $\alpha_{m_a}$ is updated using the following adaptation law:  \vspace{-0.4cm}
\begin{gather}
\label{eq:adaptive_ma}
\hat{\alpha}_{m_a}(i+1)=\hat{\alpha}_{m_a}(i)+\frac{T(m_a(i)-m_{a,d}(i))}{\kappa_4}
\end{gather}
After incorporating the predicted ADC imprecisions ($\mu_{\dot{m}_{ai}}$) into the air mass controller, the controller input is:  \vspace{-0.4cm}
\begin{gather}\label{eq:Engine_DSMC_Final_4}
\dot{m}_{ai}(i)=\frac{1}{T}[(\dot{m}+\hat{\alpha}_{m_a}(i))_{ao}(i)T+(\rho_4-1)s_4(i)\\+m_{a,d}(i+1)
-m_{a,d}(i)]-\mu_{\dot{m}_{ai}}(i)s_4(i)~~~~~~~~~~~~~~~~~~~ \nonumber
\end{gather}
In Eq.~(\ref{eq:Engine_DSMC_Final_1}), (\ref{eq:Engine_DSMC_Final_2}), (\ref{eq:Engine_DSMC_Final_3}), and (\ref{eq:Engine_DSMC_Final_4}),  $\mu_{\Delta}$, $\mu_{\dot{m}_{fc}}$, $\mu_{m_{a,d}}$, and $\mu_{\dot{m}_{ai}}$ are the estimations of propagated ADC uncertainties on control signals which are computed according to the approach explained in Fig.~\ref{fig:dsmc_with_PU}. Estimation of propagated ADC uncertainties on control signals requires the knowledge of ADC uncertainty on measured signals ($\mu_{T_{exh}}$, $\mu_{\dot{m}_f}$, $\mu_{\omega_e}$, and $\mu_{m_a}$) which are predicted using Eq.~(\ref{eq:StageIII_7}). 

\begin{figure}[h!]
\begin{center}
\includegraphics[angle=0,width= \columnwidth]{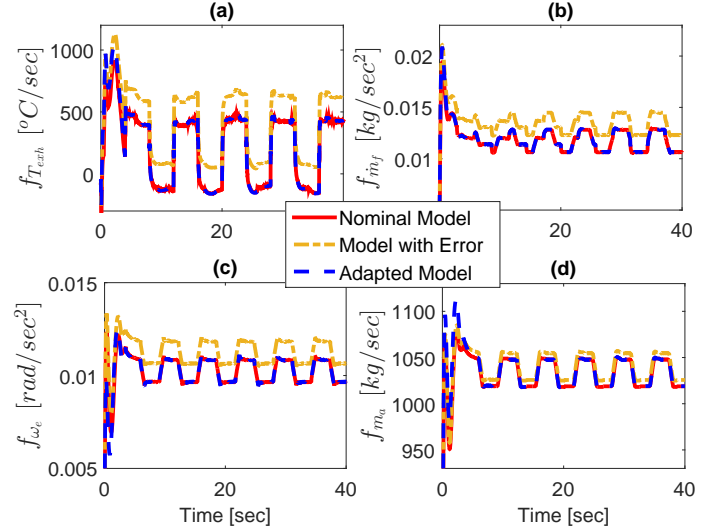} \vspace{-0.7cm}
\caption{\label{fig:Dynamics_additive_offset} The impacts of additive uncertainty term on the engine dynamics inside the DSMC and how the adaptation mechanism drives the model with error to its nominal value: (a) $T_{exh}$, (b) $\dot{m}_f$, (c) $\omega_e$, and (d) $m_a$ (20~$ms$ sampling time and 10~$bit$ quantization level).
} \vspace{-1.90cm}
\end{center}
\end{figure}

\begin{figure}[h!]
\begin{center}
\includegraphics[angle=0,width= \columnwidth]{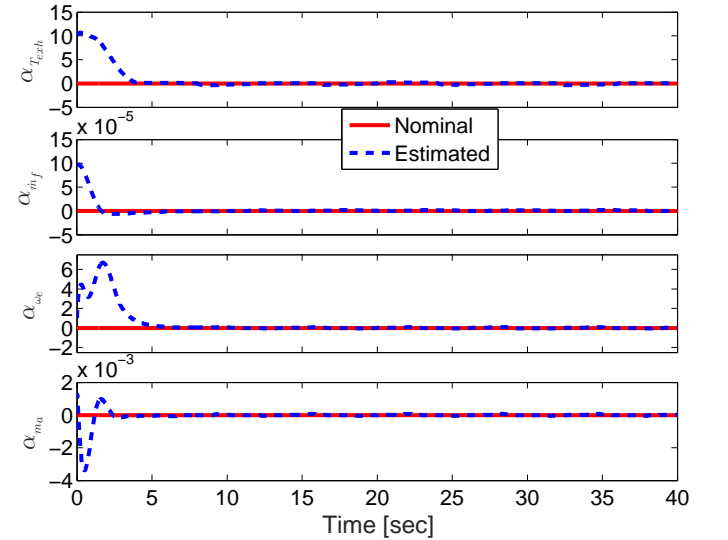} \vspace{-0.75cm}
\caption{\label{fig:Param_Conv_additive}Results of unknown additive parameters convergences (20~$ms$ sampling time and 10~$bit$ quantization level).
} \vspace{-0.95cm}
\end{center}
\end{figure}

Fig.~\ref{fig:Dynamics_additive_offset} shows the impact of unknown additive model uncertainty on part of baseline engine DSMC which represents the plant's dynamics ($f$). It is observed that in the presence of uncertainties on the model parameters, there exists a permanent error in the estimated dynamics compared to the nominal model (model with no uncertainty). This error adversely affects the tracking performance of the designed DSMC to the point that the outcome of the non-adaptive DSMC is not acceptable anymore. Once the adaptation mechanism (Eq.~(\ref{eq:Model_un_mamini14})) is activated, it can be seen from Fig.~\ref{fig:Dynamics_additive_offset} that the model with error is driven towards the nominal model and the additive model uncertainty is compensated. Fig.~\ref{fig:Param_Conv_additive} shows the results of unknown additive uncertainty term ($\hat{\alpha}$) estimation against the actual (nominal) values ($\alpha$). The adaptation mechanism removed the error in the modeled dynamics in the DSMC by driving the unknown uncertainty to its nominal value ``0" in less than 5 $sec$. 

\begin{figure}[h!]
\begin{center}
\includegraphics[angle=0,width= \columnwidth]{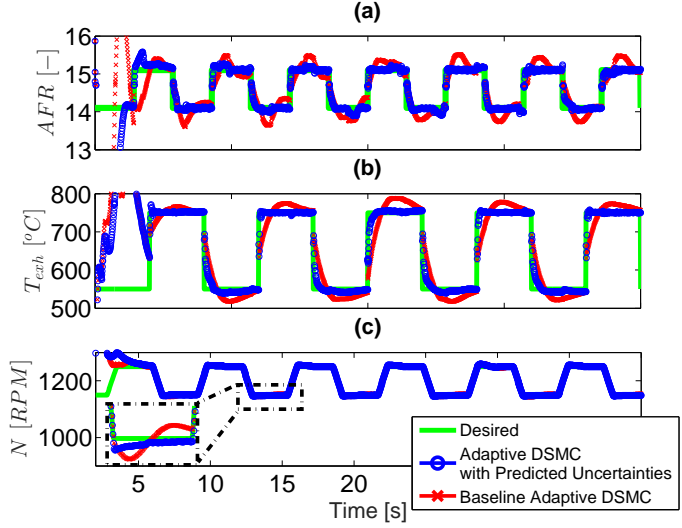} \vspace{-0.6cm}
\caption{\label{fig:EngineDSMC}Results of engine control under ADC imprecisions and additive type of model uncertainties: (a) $AFR$, (b) $T_{exh}$, and (c) $\omega_e$. (20~$ms$ sampling time and 10~$bit$ quantization level)} \vspace{-1.15cm}
\end{center}
\end{figure}

In the next stage, the predicted implementation imprecisions are incorporated in the adaptive DSMC to compare the performance of the baseline adaptive controller against adaptive DSMC with predicted ADC uncertainty. Fig.~\ref{fig:EngineDSMC} illustrates this comparison for $AFR$, $T_{exh}$, and $\omega_e$ controller, respectively. It can be clearly observed from Fig.~\ref{fig:EngineDSMC} that once the adaptation mechanism removed the additive model uncertainty within the controller (i.e. after $t$=5 $sec$), the adaptive DSMC with incorporated implementation uncertainties improves the tracking performance of the controller.\vspace{-0.85cm}
%

%$AFR$ is an important trajectory, as its accurate tracking results in meeting the desired emission and fuel consumption targets. 
When the signals at controller I/O are sampled every 20 $ms$ and the feedback signals from the plant are quantized through the ADC with 10 $bit$ quantization level, the baseline adaptive DSMC shows weak performance in tracking the desired AFR trajectory. Specifically, when the desired trajectory experiences a rapid rise or drop, the baseline DSMC cannot handle these changes properly. In addition to the AFR desired trajectory variations, due to strong correlation between the engine's state equations, any rapid changes in the engine speed profile diverges AFR from its desired value. Table~\ref{table:tracking_Results} summarizes the mean and standard deviations of the tracking errors from the baseline adaptive DSMC and adaptive DSMC with predicted ADC uncertainties. It can be seen that the adaptive DSMC with predicted ADC uncertainties improves the AFR tracking performance by 50$\%$ (in terms of mean error) compared to the baseline controller. The same improvement trends can be concluded from $T_{exh}$ and $\omega_e$ controllers in which the tracking errors are decreased by 61$\%$ and 49$\%$ upon incorporating the ADC uncertainties, respectively.   
\vspace{-0.85cm}
%\vspace{-0.5cm}
%
%
\linespread{0.8}
\begin{table*} [htbp!]
\begin{center}
\caption{Mean ($\bar{e}$) and Standard Deviation ($\sigma_e$) of Tracking Errors. Values Inside the Parentheses Show the Resulting Improvement Compared to the Baseline Adaptive DSMC.
\label{table:tracking_Results}}\vspace{0.15cm}
\begin{tabular}{lccccc}
        \hline\hline
\multicolumn{1}{c}{} & \multicolumn{2}{c}{$\bar{e}$}   &  & \multicolumn{2}{c}{$\sigma_e$} \\
        \cline{2-3} \cline{5-6}
         \textbf{}   & \textbf{Baseline }      &\textbf{Adaptive DSMC}      &  &  \textbf{Baseline }      &\textbf{Adaptive DSMC}  \\
         \textbf{}   & \textbf{Adaptive DSMC}  &\textbf{with Predicted}     &  &  \textbf{Adaptive DSMC}  &\textbf{with Predicted} \\
         \textbf{}   & \textbf{\textcolor{blue}{Reference}}&\textbf{Implementation Uncertainties} &  &  \textbf{}           &\textbf{Implementation Uncertainties}\\ \hline
             AFR [-] &      0.24            &       0.12~\textcolor{blue}{(-50.0\%)}           &  &        0.29           & 0.21   \\ \hline \vspace{0.05cm}
%                     & & \textcolor{blue}{(-50.0\%)} &  &  &  \\ \hline \vspace{0.05cm}
  $T_{exh}$~[$^o$C]   &  27.9 &  10.9~\textcolor{blue}{(-60.9\%)}                       &   & 40.0 & 29.8    \\ \hline%\vspace{0.05cm}                      %                     &       & \textcolor{blue}{(-60.9\%)} &   &      &         \\ \hline %\vspace{0.05cm}
  $N$ [RPM]          & 3.7   & 1.9~\textcolor{blue}{(-48.6\%)}                         &   & 4.2 & 2.5    \\ %\vspace{0.05cm}
 %                    &       & \textcolor{blue}{(-48.6\%)} &   &      &         \\
    \hline\hline
    \end{tabular}
\end{center}
\end{table*}
\linespread{1} \vspace{-0.5cm}
\section{Selection of Model Uncertainty Type}
\label{sec:CombinedAdaptiveDSMC}
In the previous section, it was shown that by using the proposed adaptation law in Eq.~(\ref{eq:Model_un_mamini14}), the additive type of model uncertainty can be compensated and the model with error converges to the nominal model. The adaptation theory for overcoming the multiplicative type of model uncertainty was studied in section~\ref{multiplicative_adaptive_section}. The challenging question lies in whether the error in the modeled dynamics should be considered as additive or multiplicative uncertainty. The answer to this question can be traced in the state equations of the plant model. 

In the engine control problem, the nonlinear model equations ($f_{exh},~f_{\dot{m}_f},~f_{\omega_e}$, and $f_{m_a}$) are divided into two groups: (i) rotational and intake air manifold dynamics, and (ii) exhaust gas temperature and fuel flow dynamics. As discussed in the previous section, the rotational and the intake air mass flow dynamics are obtained using torque map and volumetric efficiency curve fit, respectively. According to Eq.~(\ref{eq:Engine_discretized_we}), $\alpha_3$ indicates the error in reading the torque map ($T_{loss}$). Since reading the values from the torque map between break points in the look-up table uses a linear relationship, the additive uncertainty term can express the potential error in the modeled dynamics accordingly. The same conclusion is valid for Eq.~(\ref{eq:Engine_discretized_ma}), because the volumetric efficiency curve fit is calculated using a second order curve fit (Eq.~(\ref{eq:Engine_discretized_ma2}) and (\ref{eq:Engine_discretized_ma3})). Thus, selecting the additive type of model uncertainty is more consistent with the state equations of the rotational and air mass flow dynamics.%, and the same adaptation laws in Eq.(\ref{eq:adaptive_we}) and (\ref{eq:adaptive_ma}) are used to compensate for the error in the rotational and air mass flow dynamics, respectively.

\begin{figure}[h!]
\begin{center}
\includegraphics[angle=0,width= \columnwidth]{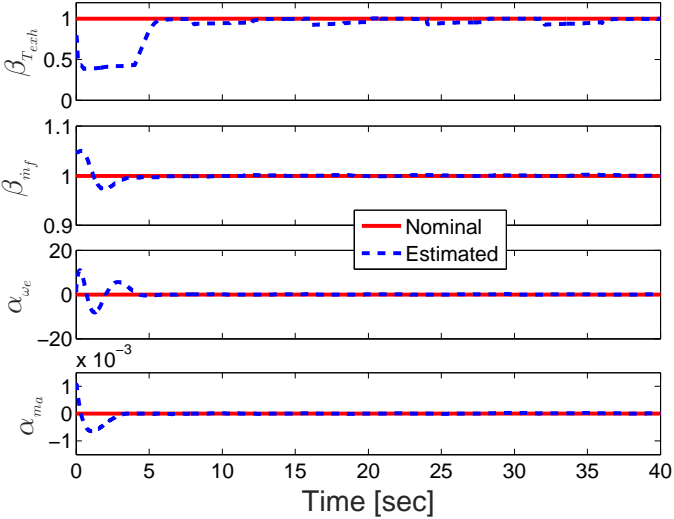} \vspace{-0.75cm}
\caption{\label{fig:Param_Conv_combined}Convergence results of unknown additive and multiplicative parameters (combined case).} \vspace{-0.7cm}
\end{center}
\end{figure}

On the other side, exhaust temperature dynamics (Eq.~(\ref{eq:Engine_discretized_Texh})) depends on the  exhaust gas time constant ($\tau_{e}$) heavily. Similarly, fuel evaporation time constant $\tau_{f}$ dictates the dynamics of the fuel flow rate into the cylinder. Thus, any error in estimating the time constants ($\tau_{e}$ and $\tau_{f}$) results in significant deviation from the nominal models. To this end, instead of assuming additive uncertainty term on $T_{exh}$ and $\dot{m}_f$ dynamics (Eq.~(\ref{eq:Engine_discretized_Texh}) and (\ref{eq:Engine_discretized_mdotf})), multiplicative uncertainty terms are assumed to compensate for any error in estimating $\tau_{e}$ and $\tau_{f}$, respectively. For $T_{exh}$ dynamics, Eq.~(\ref{eq:Engine_discretized_Texh}) is updated with the multiplicative uncertainty term ($\beta_{T_{exh}}$) as follows:  \vspace{-0.5cm}
\begin{gather}
\label{eq:Engine_discretized_Texh_multi}
f_{{T_{exh}}}=\beta_{T_{exh}}\times\frac{1}{\tau_e}[600AFI-T_{exh}]
\end{gather}
The adaptation law to drive $\beta_{T_{exh}}$ to its nominal value is obtained according to Eq.~(\ref{eq:Model_un_selina1}): \vspace{-0.4cm}
\begin{gather}
\label{eq:adaptive_Texh_multi}
\hat{\beta}_{T_{exh}}(i+1)=\hat{\beta}_{T_{exh}}(i)+\frac{T(s_1(i))}{\tau_e \rho_{\beta_1}}(600AFI-T_{exh}(i))
\end{gather}
and the new control input is: \vspace{-0.4cm}
\begin{gather}\label{eq:Engine_DSMC_Final_1_multi}
\Delta(i)=\frac{\tau_e}{7.5AFI\,.\,T}[-\hat{\beta}_{T_{exh}}(i)\frac{T}{\tau_e}(600AFI-T_{exh}(i))~~~~~~~~  \\
+(\rho_1-1)s_1(i)+T_{exh,d}(i+1)-T_{exh,d}(i)]-\mu_{\Delta}(i)s_1(i) \nonumber
\end{gather} %\vspace{-1cm}
For fuel flow dynamics, similar to $T_{exh}$ controller, Eq.~(\ref{eq:Engine_discretized_mdotf}) is updated with multiplicative uncertainty term ($\beta_{\dot{m}_f}$) to represent the error in estimating the fuel evaporation time constant ($\tau_f$): \vspace{-0.4cm}
\begin{gather}
\label{eq:Engine_discretized_mdotf_multi}
f_{\dot{m}_f}=-\beta_{\dot{m}_f}\times\frac{1}{\tau_f}(\dot{m}_f(i)) 
\end{gather} 
The adaptation becomes: \vspace{-0.4cm}
\begin{gather}
\label{eq:adaptive_mdotf_multi}
\hat{\beta}_{\dot{m}_f}(i+1)=\hat{\beta}_{\dot{m}_f}(i)+\frac{T(s_2(i))}{\tau_f \rho_{\beta_2}}\dot{m}_f(i)
\end{gather}
Finally, the control input can be expressed as: \vspace{-0.4cm}
\begin{gather}\label{eq:Engine_DSMC_Final_2_multi}
\dot{m}_{fc}(i)=\frac{\tau_f}{T}[\hat{\beta}_{\dot{m}_f}(i)\frac{T}{\tau_f}(\dot{m}_f(i))+(\rho_2-1)s_2(i)~~~~~~~~~~~~~~~  \\
+\dot{m}_{f,d}(i+1)-\dot{m}_{f,d}(i)]- \mu_{\dot{m}_{fc}}(i)s_2(i) \nonumber
\end{gather} %\vspace{-1cm}
Fiq.~\ref{fig:Param_Conv_combined} shows the unknown parameter convergences for the combined case in which uncertainties on $T_{exh}$ and $\dot{m}_f$ are multiplicative and driven to ``1", and additive uncertainty terms on $m_a$ and $\omega_e$ dynamics are driven to the nominal value ``0". Once the model uncertainty is compensated using the proposed adaptation mechanisms, the baseline adaptive DSMC can be improved by incorporating the knowledge of propagated ADC uncertainties on control signals, similar to section~\ref{subsec:Engine_Control}. In the next section, the real-time performance of the proposed adaptive DSMC with ADC uncertainties and combined adaptation laws are demonstrated on a processor-in-the-loop (PIL) setup. \vspace{-0.5cm}

\section{DSMC Real-Time Verification} \label{sec:Results}
\begin{figure}[h!]
\begin{center}
\includegraphics[angle=0,width= \columnwidth]{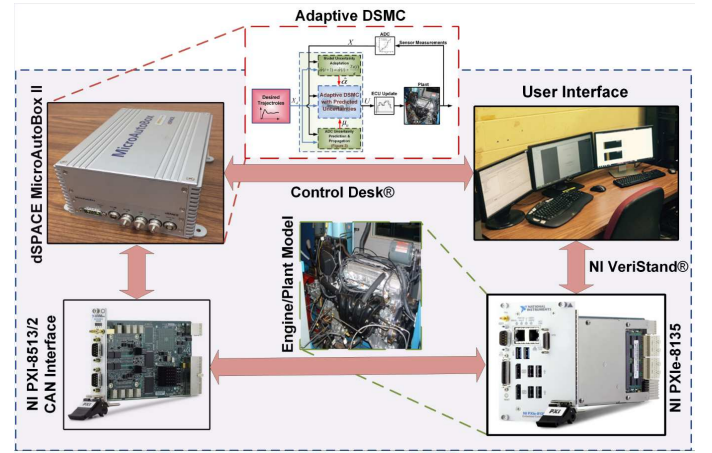} \vspace{-0.6cm}
\caption{\label{fig:HIL_PXI_Schematic}Schematic of the PIL setup for real-time DSMC verification.} \vspace{-1.75cm}
\end{center}
\end{figure}

To verify the proposed adaptive DSMC performance in real-time, the designed controller in the previous section with combined adaptation mechanisms are tested in a PIL setup shown in Fig.~\ref{fig:HIL_PXI_Schematic}. Two processors are utilized in the PIL setup. The first one is a National Instrument (NI) PXI processor (NI PXIe-8135) in which the generated C-code of the engine plant model is implemented. The second processor is a dSPACE MicroAutoboxII (MABX) which represents the main engine controller unit (ECU). The generated C-code of the adaptive DSMC along with the adaptation and uncertainty prediction mechanisms are implemented into MABX. The trajectory tracking performance of the controller and the engine operation are tested real-time under embedded ADC imprecisions (sampling time=20 $ms$ and quantization level=10 $bit$) on feedback and control signals, in the presence of additive and multiplicative types of model uncertainties in the controller. %The reason for studying the controller performance under 10 $ms$ of sampling time and quantization level of 10 $bit$ is that these values were recommended as the ''minimally implemented" controller in~\cite{Amini_DSC,AminiSAE2016}.
Real-time test configuration is conducted using NI VeriStand$^{\textregistered}$~and dSPACE Control Desk$^{\textregistered}$~software on an interface desktop computer. %\vspace{-0.75cm}

\begin{figure}[h!]
\begin{center}
\includegraphics[angle=0,width= \columnwidth]{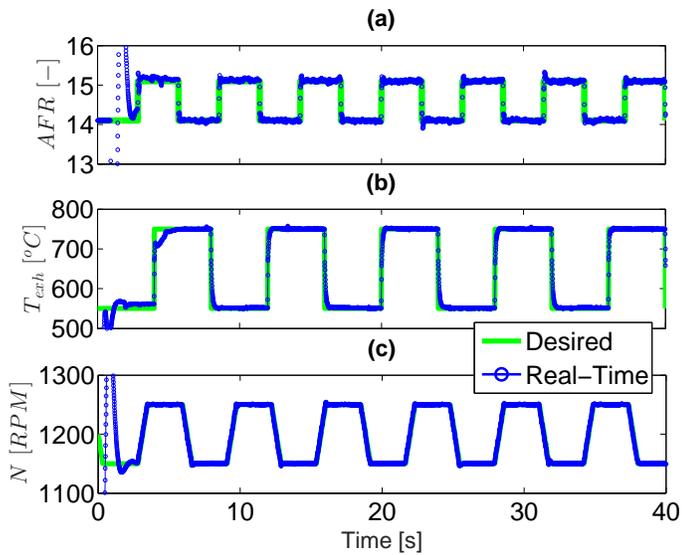} \vspace{-0.5cm}
\caption{\label{fig:PIL_Test}Real-time verification for the engine control using the proposed
DMSC under 20 $ms$ sampling time and 10 $bit$ quantization level (Additive type of model uncertainty is used for rotational and intake air mass flow dynamics, and multiplicative type is used for fuel flow and exhaust gas temperature dynamics).} \vspace{-0.75cm}
\end{center}
\end{figure}

Fig.~\ref{fig:PIL_Test} shows the results of real-time PIL testing of the proposed adaptive DSMC with incorporated ADC uncertainties on control signals for the engine case study. The results verify that the adaptive DSMC is able to (i) remove the uncertainty in the model in less than 5 seconds, (ii) track all the desired trajectories under ADC uncertainties, and (iii) operate in real-time, since it is computationally efficient. \vspace{-0.75cm}
\section{Summary and Conclusions}   \label{sec:Conclusion}
In this paper, we presented a generic discrete sliding mode controller design methodology with consideration of implementation (ADC) imprecisions and model uncertainty. First, an online methodology was developed for predicting and propagating the sampling and quantization imprecisions on both measured and control signals. In the next step, the DSMC was examined under two different types of model uncertainties (multiplicative and additive). A discrete Lyapunov argument was employed to evaluate the DSMC stability with model uncertainty and also, to derive the adaptation laws to compensate the unknown additive uncertainty term in the controller equations. The final controller is an adaptive DSMC with incorporated ADC uncertainties which can handle model and implementation uncertainties online. A real-time PIL testing was conducted to verify the performance of the proposed adaptive DSMC on an automotive engine tracking control problem. The PIL testing on the actual ECU showed that the new adaptive DSMC design is able to remove the error in the modeled dynamics quickly. Additionally, upon incorporation of ADC knowledge into the adaptive DSMC, the tracking performance of the controller under sampling and quantization imprecisions improved by 50-60\% compared to a baseline DSMC design.\vspace{-0.5cm}
%
%%%%%%%%%%%%%%%%%%%%%%%%%%%%%%%%%%%%%%%%%%%%%%%%%%%%%%%%%%%%%%%%%%%%%%%
\begin{acknowledgment}
This material is based upon the work supported by the National Science Foundation under Grant No. 1434273. Dr. Ken Butts from Toyota Motor Engineering $\&$ Manufacturing North America is gratefully acknowledged for his technical comments during the course of this study.
\end{acknowledgment} \vspace{-0.85cm}
%%%%%%%%%%%%%%%%%%%%%%%%%%%%%%%%%%%%%%%%%%%%%%%%%%%%%%%%%%%%%%%%%%%%%%
% The bibliography is stored in an external database file
% in the BibTeX format (file_name.bib).  The bibliography is
% created by the following command and it will appear in this
% position in the document. You may, of course, create your
% own bibliography by using thebibliography environment as in
%
% \begin{thebibliography}{12}
% ...
% \bibitem{itemreference} D. E. Knudsen.
% {\em 1966 World Bnus Almanac.}
% {Permafrost Press, Novosibirsk.}
% ...
% \end{thebibliography}

% Here's where you specify the bibliography database file.
% The full file name of the bibliography database for this
% article is asme2e.bib. The name for your database is up
% to you.
%\bibliography{asme2e}
\bibliographystyle{IEEEtran}
\bibliography{DSMC_UP_bib}

%%%%%%%%%%%%%%%%%%%%%%%%%%%%%%%%%%%%%%%%%%%%%%%%%%%%%%%%%%%%%%%%%%%%%%
%\appendix       %%% starting appendix
%\section*{Appendix A: Head of First Appendix}
%Avoid Appendices if possible.
%
%%%%%%%%%%%%%%%%%%%%%%%%%%%%%%%%%%%%%%%%%%%%%%%%%%%%%%%%%%%%%%%%%%%%%%%
%\section*{Appendix B: Head of Second Appendix}
%\subsection*{Subsection head in appendix}
%The equation counter is not reset in an appendix and the
%numbers will follow one continual sequence from the
%beginning of the article to the very end as shown in the
%following example.
%\begin{equation}
%a = b + c.
%\end{equation}

\end{document}